\documentclass{article}\begin{document}
\centerline{\bf On the solution of the general quintic}
\centerline{\bf using the Rogers-Ramanujan continued fraction}
\vskip .4in

\centerline{Nikos Bagis}

\centerline{Stenimahou 5 Edessa}
\centerline{Pella 58200, Greece}
\centerline{bagkis@hotmail.com}
\vskip .2in

\[
\]


\[
\]

\centerline{\bf Abstract}

\begin{quote}
In this article we give solution of the general quintic equation by means of the Rogers-Ramanujan continued fraction. More precisely we express a root of the quintic as a known algebraic function of the Rogers-Ramanujan continued fraction.     
\end{quote}

\section{Introduction and Definitions}

In the most general case the quintic equation read as 
\begin{equation}
ax^5+bx^4+cx^3+dx^2+ex+f=0\textrm{, }a\neq0
\end{equation}
According to Abel's impossibility theorem (see [22]), this equation has no general solution. With the word solution we mean a known algebraic function of its coefficients. For example the quadratic equation
\begin{equation}
ax^2+bx+c=0\textrm{, }a\neq0
\end{equation} 
have two solutions
\begin{equation}
x_{1,2}=\frac{-b\pm\sqrt{b^2-4ac}}{2a}
\end{equation}
The cubic equation 
\begin{equation}
ax^3+bx^2+cx+d=0\textrm{, }a\neq0
\end{equation}
admits three solutions, one necessary real and two other complex in general (see Appendix for a very amusing way of solving (4) with the help of (2)).\\The quartic equation 
\begin{equation}
ax^4+bx^3+cx^2+dx+e=0\textrm{, }a\neq0
\end{equation}
has solution also. For higher degrees, greater than four we can not construct a general solution.\\
In theory of elliptic functions (see [20]) exist functions that take algebraic values. These functions are solutions of algebraic polynomial equations that we can not solve. For example, the null Jacobi theta functions
\begin{equation}
\vartheta_2(q):=\sum^{\infty}_{n=-\infty}q^{(n+1/2)^2}\textrm{, }\vartheta_3(q):=\sum^{\infty}_{n=-\infty}q^{n^2}\textrm{, }\vartheta_4(q):=\sum^{\infty}_{n=-\infty}(-1)^nq^{n^2},
\end{equation}
where $|q|<1$, are forming the elliptic singular modulus 
\begin{equation}
k_r:=\frac{\vartheta_2^2(q)}{\vartheta_3^2(q)}\textrm{, }q=e^{-\pi\sqrt{r}}\textrm{, }r>0.
\end{equation}
The elliptic singular modulus take algebraic values for all $r-$positive rationals. A few values of $k_r$ are $k_1=\frac{1}{\sqrt{2}}$, $k_2=\sqrt{2}-1$, $k_3=\frac{1}{2}\sqrt{2-\sqrt{3}}$, $k_4=3-2\sqrt{2}$, $k_5=\sqrt{\frac{9+4\sqrt{5}-2\sqrt{38+17\sqrt{5}}}{18+8\sqrt{5}}}$, $k_{2/3}=-3+2\sqrt{3}+\sqrt{14-8\sqrt{3}}$, $\ldots$. The complexity of these algebraic numbers is getting bigger as we move to greater and more involved $r$. An example is $k_{6/7}$, which is root of 
$$
1-6600x+51988x^2-6600x^3-103926x^4+6600x^5+
$$
\begin{equation}
+51988x^6+6600x^7+x^8=0.
\end{equation} 
This happens to be symmetrical octic, hence it is solvable.\\ 
For a list of values of $k_r$ one can see [21],[15],[18] and Wolfram pages in the Web. We may also, knowing the value of $k_r$ proceed to higher values using solvable modular equations. For example hold the following relations ([7] chapter 18, pg.213-219)
\begin{equation}     
k_{4r}=\frac{1-\sqrt{1-k_r^2}}{1+\sqrt{1-k_r^2}}\textrm{, }\forall r>0
\end{equation}
and ([7] chapter 19, Entry 5, pg.230)
\begin{equation}
\sqrt{k_rk_{9r}}+\sqrt{k'_rk'_{9r}}=1,
\end{equation}
where $k'_r=\sqrt{1-k_r^2}$, $r>0$, ((9) and (10) are the solvable modular equations of degree's 2 and 3 respectively of the singular modulus).\\
The elliptic singular modulus $k_r$ is also the solution $x$ of the equation
\begin{equation}
\frac{{}_2F_{1}\left(\frac{1}{2},\frac{1}{2};1;1-x^2\right)}{{}_2F_{1}\left(\frac{1}{2},\frac{1}{2};1;x^2\right)}
=\sqrt{r},
\end{equation}
where ${}_2F_1\left(a,b;c;x\right)$ is the known Gauss hypergeometric function (see [23]). In particular
\begin{equation}
{}_2F_{1}\left(\frac{1}{2},\frac{1}{2};1;x^2\right)=\sum^{\infty}_{n=0}\frac{\left(\frac{1}{2}\right)^2_n}{(n!)^2} x^{2n}=\frac{2}{\pi}K(x)=\frac{2}{\pi}\int^{\pi/2}_{0}\frac{d\phi}{\sqrt{1-x^2\sin^2(\phi)}}
\end{equation}

The 5th degree modular equation which connects $k_{25r}$ and $k_r$, (see [7]) is:
\begin{equation}
k_rk_{25r}+k'_rk'_{25r}+2^{5/3} (k_rk_{25r}k'_rk'_{25r})^{1/3}=1
\end{equation} 
The problem of solving (13) and finding $k_{25r}$ reduces to that of solving what Hermite called depressed equation (see [20]):
\begin{equation}
u^6-v^6+5u^2v^2(u^2-v^2)+4uv(1-u^4v^4)=0,
\end{equation}
where $u=k^{1/4}_{25r}$ and $v=k^{1/4}_{r}$.\\
The depressed equation is a sextic equation, which in general is not solvable in radicals and also very important for the solution of the quintic (1).  
With the help of null theta functions we can obtain a closed form solution
\begin{equation}
k_{25r}=\frac{\vartheta^2_2(q^5)}{\vartheta^2_3(q^5)}.
\end{equation} 
But this is not satisfactory, since finding values of the $\vartheta$ functions, is the same as using the $k_r$ and (14) (we have no further representations of $\vartheta$ functions to lead us in some general theory). That is, we again need to know from the value of $k_r$, the value of $k_{25r}$.\\
In this article we give a parametrization for the pairs of solutions $\left\{k_{r},k_{25r}\right\}$ (see (47),(48) below) in the form $\left\{A_1(R(q)),A_2(R(q))\right\}$, $q=e^{-\pi\sqrt{r}}$, with $A_1(x)$ and $A_2(x)$ known algebraic functions and where $R(q)$ is the Rogers-Ramanujan continued fraction (RRCF)
\begin{equation}
R(q)=\frac{q^{1/5}}{1+}\frac{q^1}{1+}\frac{q^2}{1+}\frac{q^3}{1+}\ldots\textrm{, }|q|<1.
\end{equation} 
In the question ''Why the singular modulus must be replaced by the (RRCF) ?'' the answer is:\\   
As like the singular modulus, the (RRCF) obeys modular equations. The 2nd, 3rd and 5th degree modular equations of RRCF are (see [19],[11]):
\begin{equation}
\frac{R\left(q^2\right)-R\left(q\right)^2}{R\left(q^2\right)+R\left(q\right)^2}=R(q)R\left(q^2\right)^2,
\end{equation}
\begin{equation}
\left(R\left(q^3\right)-R\left(q\right)^3\right)\left(1+R(q)R\left(q^3\right)^3\right)=3R(q)^2R\left(q^3\right)^2,
\end{equation}
\begin{equation}
R\left(q^{1/5}\right)^5=R(q)\frac{1-2R(q)+4R(q)^2-3R(q)^3+R(q)^4}{1+3R(q)+4R(q)^2+2R(q)^3+R(q)^4}
\end{equation}
Observe that the 5th degree modular equation is solved already and we have a very strong advantage against the elliptic singular modulus $k_r$. Moreover as we will show knowing for a certain $r$ the value of $R(q)$, $q=e^{-\pi\sqrt{r}}$, $r>0$ we know all $k_{25^nr}$, $n\in \textbf{N}$ (although in view of [25] we can extend $\textbf{N}$ into $\textbf{Z}$).

\section{Main Results}

We begin constructing the functions $A_1(x)$ and $A_2(x)$.\\ 
\\
We call $T_1(x)$ the algebraic function that relates $R\left(q^2\right)$ and $R(q)$. Set $T_4(x)$ to be the algebraic function that relates $R\left(q^{1/5}\right)$ and $R(q)$. Hence
\begin{equation}
R\left(q^2\right)=T_1\left(R(q)\right)\textrm{, }R\left(q^{1/5}\right)=T_4\left(R(q)\right) 
\end{equation}
If we know the value of (RRCF), then we can find the value of $j-$invariant from Klein's equation (see [24] and Wolfram pages 'Rogers Ramanujan Continued fraction'):
\begin{equation}
j_r=j(q)=-\frac{\left(R_2^{20}-228 R_2^{15}+494 R_2^{10}+228 R_2^5+1\right)^3}{R_2^5 \left(R_2^{10}+11 R_2^5-1\right)^5},
\end{equation}
where $R_2=R(q^2)=T_1\left(R(q)\right)$. Hence we set $T_2(x)$ such
\begin{equation}
j_r=j(q)=T_2\left(R(q^2)\right)=T_2\circ T_1\circ R(q)
\end{equation}
From the relation 
\begin{equation}
j_r=\frac{256(k^2_r+(k'_{r})^4)^3}{(k_rk'_r)^4} , 
\end{equation}
we can set $T_3(x)$ such that
\begin{equation}
k_r=T_3\left(j(q)\right)=T_3\circ T_2\circ T_1\circ R(q).
\end{equation}
Note here that (23) is solvable with respect to $k_r$ in radicals. In the method described above one must be very careful which branch of the solution (of the polynomial equation) is using.\\

Knowing $k_r$ and $k_{r/25}$, we can evaluate $k_{25r}$ (see [25]) and specify relations of the form
\begin{equation}
k_{25r}=\Psi\left(k_r,k_{r/25}\right)\textrm{ and }k_{25^nr}=\Psi_{n}\left(k_r,k_{r/25}\right)\textrm{, }n\in\textbf{N}
\end{equation} 
The next theorem describes in details the relations (25).\\
\\ 
\textbf{Theorem 1.} (Bagis [25])\\
For all $n=1,2,\ldots$ we have
\begin{equation}
k_{25^nr_0}=\sqrt{1/2-1/2\sqrt{1-4\left(k_{r_0}k'_{r_0}\right)^2\prod^{n}_{j=1}P^{(j)}\left[\sqrt[12]{\frac{k_{r_0}k'_{r_0}}{k_{r_0/25}k'_{r_0/25}}}\right]^{24}}},
\end{equation}
where the function $P$ is in radicals known function and is given from
\begin{equation}
P(x)=P[x]=U\left[Q\left[U^{*}\left[x\right]^6\right]^{1/6}\right]
\end{equation}
\begin{equation}
P^{(n)}(x)=(P\underbrace{\circ\ldots\circ}_{n} P)(x)\textrm{ , } P^{(1)}(x)=P(x) . 
\end{equation}
The function $Q$ is that of (29) and $U$, $U^{*}$ are that of (32) and (33) below.\\      
\\

We give the analytic expansion of $Q$: 
\begin{equation}
Q(x)=\frac{\left(-1-e^{\frac{1}{5}y}+e^{\frac{2}{5} y}\right)^5}{\left(e^{\frac{1}{5}y}-e^{\frac{2}{5}y}+2 e^{\frac{3}{5} y}-3 e^{\frac{4}{5}y}+5 e^{y}+3 e^{\frac{6}{5}y}+2 e^{\frac{7}{5}y}+e^{\frac{8}{5} y}+e^{\frac{9}{5}y}\right)}
\end{equation}
and
\begin{equation}
y=\textrm{arcsinh}\left(\frac{11+x}{2}\right).
\end{equation}
Note that (for no confusion), the function $\exp\left(rational\cdot y\right)$ is algebraic function of $x$, since $y$ is the arc sine of a rational function.\\   
Consider now the following equation
\begin{equation}
\frac{X^2}{\sqrt{5}Y}-\frac{\sqrt{5}Y}{X^2}=\frac{1}{\sqrt{5}}\left(Y^3-Y^{-3}\right)
\end{equation}
This equation is solvable in radicals with respect to $Y$ and $X$ also. One can find the solution
\begin{equation}
Y=U(X)=\sqrt{-\frac{5}{3 X^2}+\frac{25}{3 X^2 h(X)}+\frac{X^4}{h(X)}+\frac{h(X)}{3 X^2}}
\end{equation}
where
$$
h(x)=\left(-125-9x^6+3\sqrt{3}\sqrt{-125x^6-22x^{12}-x^{18}}\right)^{1/3}
$$
The solution of (31) with respect to $X$ is
\begin{equation}
X=U^{*}(Y)=\sqrt{-\frac{1}{2 Y^2}+\frac{Y^4}{2}+\frac{\sqrt{1+18 Y^6+Y^{12}}}{2 Y^2}} .
\end{equation}
\\

We can simplify Theorem 1 replacing it with an equivalent transformation.\\
Set $p$ such that
\begin{equation}
U^{*}(x)^6=-11+2\sinh\left(5\log\left(\frac{p+\sqrt{4+p^2}}{2}\right)\right),
\end{equation}
with
\begin{equation}
p=2\sinh\left(\frac{1}{5}\textrm{arcsinh}\left(\frac{11+U^{*}(x)^6}{2}\right)\right).
\end{equation}
Then (29) becomes
\begin{equation}
Q(p)=\frac{(p-1)^5}{11+6p+6p^2+p^3+p^4}.
\end{equation}
Hence if 
\begin{equation}
s=\sqrt[3]{\frac{(p-1)^5}{11+6p+6t^2+p^3+p^4}},
\end{equation}
then $P_1=P(x)=\sqrt{Y}$, where $Y$ is root of 
\begin{equation}
Y^3+5s^{-1}Y^2-s Y-1=0.
\end{equation}
\\
Hence we get the next:\\
\\ 
\textbf{Lemma 1.}\\
Given $k_{r_0}$ and $k_{r_0/25}$, if 
\begin{equation}
\alpha=U^{*}\left(\sqrt[12]{\frac{k_{r_0}k'_{r_0}}{k_{r_0/25}k'_{r_0/25}}}\right)^6
\end{equation}
and
\begin{equation}
p=2\sinh\left(\frac{1}{5}\textrm{arcsinh}\left(\frac{11+\alpha}{2}\right)\right),
\end{equation}
then
\begin{equation}
s=\sqrt[3]{\frac{(p-1)^5}{11+6p+6t^2+p^3+p^4}}.
\end{equation}
Hence
\begin{equation}
P\left(\sqrt[12]{\frac{k_{r_0}k'_{r_0}}{k_{r_0/25}k'_{r_0/25}}}\right)=\sqrt{Y(s)},
\end{equation}
where 
\begin{equation}
Y^3+5s^{-1}Y^2-sY-1=0\textrm{, }Y=Y(s).
\end{equation}  
\\

From the above Lemma 1 we get the value of $k_{25r}$ knowing $k_r$ and $k_{r/25}$:\\
\\
\textbf{Proposition 1.}
\begin{equation}
k_{25r}=\Psi\left(k_r,k_{r/25}\right)=\sqrt{\frac{1}{2}-\frac{1}{2}\sqrt{1-4(k_{r}k'_{r})^2\cdot  Y\left(\sqrt[3]{\frac{(p-1)^5}{11+6p+6p^2+p^3+p^4}}\right)^{12}}}
\end{equation}
\\
\textbf{Example.}\\
In case of $r_0=25$, we need the value of   
$$
P_1=P\left(\sqrt[12]{\frac{k_{r_0}k'_{r_0}}{k_{r_0/25}k'_{r_0/25}}}\right)\stackrel{r_0=25}{=}P\left(\frac{\sqrt{5}-1}{2}\right).
$$
But 
$\alpha=U^{*}\left(\frac{\sqrt{5}-1}{2}\right)^6=-2+\sqrt{5}$, hence
$Q\left(U^{*}\left(\frac{\sqrt{5}-1}{2}\right)^6\right)=Q\left(-2+\sqrt{5}\right)$. 
With 
$$
x=a+b\sqrt{d}=-2+\sqrt{5}
$$
we get
$$
s^{1/2}=Q\left(-2+\sqrt{5}\right)^{1/6}=\sqrt[6]{\frac{(p-1)^5}{11+6p+6p^2+p^3+p^4}}
$$
where
$$
p=2\sinh\left(\frac{1}{5}\textrm{arcsinh}\left(\frac{9+\sqrt{5}}{2}\right)\right)
$$
Hence we take the value of $k_{625}$:
$$
k_{625}=\sqrt{\frac{1}{2}-\frac{1}{2}\sqrt{1-(51841-23184\sqrt{5})\cdot  Y\left(\sqrt[3]{\frac{(p-1)^5}{11+6p+6p^2+p^3+p^4}}\right)^{12}}}
$$
\\

For to make things more clear and simple we use the function $m(q)=m_1(z)$, $q=e^{i\pi z}$, which is related with $k_r$ by
\begin{equation}
m_1(i\sqrt{r})=k_r\textrm{, }q=e^{-\pi \sqrt{r}},
\end{equation}
when $r>0$ and for complex values is
\begin{equation}
m(z):=k_{\left(\frac{\log(q)}{\pi}\right)^2}\textrm{, }m_1'(z)=\sqrt{1-m_1(z)^2},
\end{equation}
when $q=e^{i\pi z}$, $Im(z)>0$.\\
The first equation of (46) in better detail (since $k_r$ is defined for $r$ positive real) means that $m_1(z)$ is actualy a generalization of the singular modulus and solution of the equation
\begin{equation}
i\frac{K\left(\sqrt{1-m_1(z)^2}\right)}{K\left(m_1(z)\right)}=i\frac{K'}{K}=z\textrm{, }q=e^{i\pi z}.
\end{equation}
when $Im(z)>0$ and $-1<Re(z)\leq 1$. Also then
\begin{equation}
m_1(z)=\frac{\vartheta^2_2(q)}{\vartheta^2_3(q)}\textrm{, where }q=e^{i\pi z}.
\end{equation}
\\
\textbf{Proposition 2.}\\
The solution of the depressed equation in terms of the Rogers-Ramanujan continued fraction is:
\begin{equation}
t=m_1(z)=\left(T_{36}\circ T_2\circ T_1\right)\left(v(z)\right)\textrm{, }
l=m_1\left(z/5\right)=\left(T_{38}\circ T_2\circ T_4\circ T_1\right)\left(v(z)\right)
\end{equation}
and
\begin{equation}
m_1\left(5z\right)=\Psi(t,l),
\end{equation}
where $v(z)=R(q)$, $q=e^{i\pi z}$, $Im(z)>0$.\\
\\
\textbf{Example.}\\
Suppose $r=16$ following classical theory we have 
$$
k_{16}=33+24 \sqrt{2}-4 \sqrt{140+99 \sqrt{2}}
$$
The value of $k_{25\cdot16}=k_{400}$ is symmetric and of degree 16, which is quite difficult to evaluate. But it happens that
$$
R\left(e^{-2\pi}\right)=-\frac{1+\sqrt{5}}{2}+\sqrt{\frac{5+\sqrt{5}}{2}}
$$
hence
$$
v_{16}=R\left(e^{-4\pi }\right)=T_1\left(R\left(e^{-2\pi}\right)\right)
$$ 
and from (49) the values of $t=k_{16}$ and $l=k_{16/25}$ are known. From (44) we find the value of $k_{400}$ (in radicals closed form).\\  
\\

Hermite has proved (see [20]) the next\\
\\
\textbf{Theorem 2.} (Hermite)\\
The equation
\begin{equation}
x^5-x-a=0
\end{equation}
with
\begin{equation}
a=\frac{2(1+m_1(z)^2)}{\sqrt[4]{5^5}\sqrt{m_1(z)}m_1'(z)}
\end{equation}
and $z=i\sqrt{r}$, $r>1$ have real root
\begin{equation}
x=\frac{\Phi(z)}{2\sqrt[4]{5^3}m_1(z)^{1/4}m_1'(z)},
\end{equation}
where
$$
\Phi(z)=\left[m_1\left(5z\right)^{1/4}+m_1\left(z/5\right)^{1/4}\right]\left[m_1\left(\frac{z+16}{5}\right)^{1/4}\zeta_{1/2}-m_1\left(\frac{z+64}{5}\right)^{1/4}\overline{\zeta}_{1/2}\right]\times
$$
\begin{equation}
\times\left[m_1\left(\frac{z+32}{5}\right)^{1/4}\overline{\zeta}_{1/4}-m_1\left(\frac{z+48}{5}\right)^{1/4}\zeta_{1/4}\right]
\end{equation}
and 
$$
\phi(z)=m_1(z)^{1/4}=\sqrt{\frac{\vartheta_2(q)}{\vartheta_3(q)}}\textrm{, }q=e^{\pi i z}\textrm{, }Im(z)>0. 
$$
Also $\zeta_{\nu}=e^{i\pi \nu}$ and $\overline{w}$ denotes the complex conjugate of $w$.
\\

In the above Hermite's theorem we can write $\Phi$ equivalently as (here we use $m(q)=m_1(z)$, $q=e^{i\pi z}$):
$$
\Phi(z)=4\left[m\left(q^{5}\right)^{1/4}+m\left(q^{1/5}\right)^{1/4}\right]Im\left[m\left(-q^{1/5}\zeta_{1/5}\right)^{1/4}\zeta_{1/2}\right]\times
$$
\begin{equation}
\times Re\left[m\left(q^{1/5}\zeta_{2/5}\right)^{1/4}\zeta_{1/4}\right],
\end{equation}
where $Re(z)$, $Im(z)$ are the real and imaginary part of $z$ respectively.\\ 
Therefore we have to evaluate 
\begin{equation}
m\left(-(-q)^{1/5}\right)=m\left(-(-1)^{1/5}q^{1/5}\right)
\end{equation}
and
$$
m\left((-1)^{2/5}q^{1/5}\right)=m\left((-1)^{1/5}(-q)^{1/5}\right)\eqno{(56.1)}
$$
This can be done using
\begin{equation}
m(-q)=\frac{m(q)}{\sqrt{m(q)^2-1}}.
\end{equation}
This means that we can define the function $T_5$ such that
\begin{equation}
m(-q)=T_5\left(m(q)\right).
\end{equation}
Also we have that $m\left((-1)^{2/5}q^{1/5}\right)$ is known algebraic function of 
$$
m\left[\left((-1)^{2/5}q^{1/5}\right)^5\right]=m(q).
$$ 
We will use the theory of Ramanujan Quantities (see [26]) for the evaluation of $m\left((-q)^{1/5}\right)$.\\
\\
If $v=R(1,2,5;q)=R(1,3,5;q)=R(q)$ is the Rogers Ramanujan continued fraction and $u=u(q)=R(1,3,10;q)$, then we conjecture that  
\begin{equation}
u^3-uv+u^2v^3+v^4\stackrel{?}{=}0
\end{equation}
\textbf{Note.} In general its an observation that different Ramanujan Quantities are related with polynomial relations (see [10],[26]):\\
\begin{equation}
k_r^2=16H^2\left(\frac{1-H^2}{1+H^2}\right)^2=\frac{(1-T)(3+T)^3}{(1+T)(3-T)^3} , 
\end{equation}
where $T=\sqrt{1-8V^3}$ and $V=V(q)=R(1,3,6;q)$ is the Ramanujan's cubic continued fraction. The quantity $H=H(q)=R(1,3,8;q)$ is the Ramanujan-Gollnitz-Gordon continued fraction.\\

The $u-$quantity satisfies Theorem 8 of [26], since $2|10$ and $(1,2)=1$, $(3,2)=1$, we have 
\begin{equation}
uu^{*}=u_2=u\left(q^2\right),  
\end{equation}
where $u^{*}=u(-q)$.\\
Another evaluation, using Mathematica is
$$
-v+(-1)^{1/5}v^*-(-1)^{1/5}v^5v^*+5(-1)^{2/5}v^4(v^*)^2-10 (-1)^{3/5}v^3(v^*)^3+
$$
\begin{equation}
+5(-1)^{4/5}v^2(v^*)^4+v(v^*)^5+v^6(v^*)^5-(-1)^{1/5}v^5(v^*)^6\stackrel{?}{=}0,
\end{equation}
where $v^{*}=R(-q)=v(-q)$.\\
Relation (62) follows with the same method as (59) (see Routine pages 8-9 of [26]).\\
\\
\textbf{Lemma 2.}\\ 
If $R(q)$ is the Rogers Ramanujan continued fraction, then 
\begin{equation}
R(1,3,10;q)\stackrel{?}{=}R(q)R\left(q^2\right).
\end{equation}
\textbf{Proof.}(Partial)\\ 
We can write
\begin{equation}
\frac{1}{n}\log(R(a,b,p;q^n))=Q\log(q)-\sum^{\infty}_{k=1}\frac{q^{nk}}{nk}\sum_{d|k}X(a,b,p;d)d , 
\end{equation} 
Multiplying with $\epsilon_n$ and summing we get  
\begin{equation}
\sum^{\infty}_{n=1}\frac{\epsilon(n)}{n}\log(R(a,b,p;q^n))=Q\log(q)\sum^{\infty}_{n=1}\epsilon(n)-\sum^{\infty}_{n=1}\frac{q^n}{n}\sum_{d|n}\tau(d)\epsilon(n/d).
\end{equation}
We make use of the properties of the Dirichlet multiplication (see [3]).\\ Assume that exists arithmetic function $X_1(n)n$ such that 
$$
X_1(n)n\ast 1=\tau(n)\ast \epsilon(n)=X(n)n\ast 1\ast \epsilon(n) , 
$$
then 
\begin{equation}
X_1(n)n=\epsilon(n)\ast X(n)n.
\end{equation} 
Assuming the $X(n)$ and $X_1(n)$ correspond to Ramanujan Quantities we seek $\epsilon(n)$ such that is always 0 after a few terms. An example is when $X(n)=X(1,3,5,n)$ and $X_1(n)=X(1,3,10,n)$. Then 
\begin{equation}
\epsilon(n)\stackrel{?}{=}1\textrm{, if }n=1,2\textrm{ and } \epsilon(n)\stackrel{?}{=}0\textrm{, if }n>2
\end{equation}
From the above, using (61),(65),(67), we have the partial proof of (63).\\
\\      
\textbf{Lemma 3.}\\
$i)$ The function $w=R(q)R^2(q^2)$ can always evaluated in radicals with respect to $v=R(q)$. The equation that relates the two functions is
\begin{equation}
(w-1)\sqrt{R(q)w}+R(q)^3(1+w)=0.
\end{equation}
Moreover\\
$ii)$ The function $v'=\left|R(-q)\right|$ can always evaluated in radical function of $R(q)$.\\ 
\\
\textbf{Proof.}\\ 
The proof of $(i)$ follows from  
$u^2=R(1,3,10;q)^2=R(q)R(q)R^2(q^2)$, hence $w=u^2/R(q)=u^2/v$, or $u=\sqrt{wv}$. Setting this identity to (59) we get the result.\\
The proof of $(ii)$ follows from (63) and (59): 
$$
u(q^2)=u(q)|u(-q)|=R(q)R(q^2)|R(-q)| R(q^2)=|R(-q)|w
$$
Suppose the solution of (59) is $u(q)=f_1(R(q))$, then $u(q^2)=f_1(R(q^2))=|R(-q)|w$, hence 
$$
|R(-q)|=\frac{f_1\left(R\left(q^2\right)\right)}{w}
$$
and the proof is complete.\\
\\
From the above Lemmas we get the next\\
\\
\textbf{Proposition 3.}\\
\begin{equation}
R(-q)\stackrel{?}{=}e^{i\pi/5}\frac{f_1\left(T_1\left(R\left(q\right)\right)\right)}{R(q)T_1\left(R(q)\right)^2}=T_6\left(R(q)\right)
\end{equation} 
where $f_1(x)$ is the solution $u=R(1,3,10;q)=f_1\left(R(q)\right)$ of (59).
\\

Note that the proof is not complete since we not prove (59) and (67), but we include it to show how we arrive to (69).\\
\\
\textbf{Proposition 4.}\\
Let $q=e^{-\pi\sqrt{r}}$ with $r$ real $1<r\leq 1+\frac{1}{15}$, then
\begin{equation}
m\left(-(-q)^{1/5}\right)=-T_5\circ t_0\circ T_{312}\circ T_2\circ T_4\circ T_6 \circ R(q)
\end{equation}
and
\begin{equation}
m\left((-1)^{2/5}q^{1/5}\right)=t_0\circ t_0\circ T_{38}\circ T_2\circ T_4\circ T_6\circ T_{7}\circ R(q),
\end{equation}
where 
$$
t_0(x)=\frac{1-\sqrt{1-x^2}}{1+\sqrt{1-x^2}}.
$$
\\
\textbf{Main Theorem.}\\
If $1<r\leq1+1/15$, $q=e^{-\pi\sqrt{r}}$, evaluate $t,l$ from (72) and (73) below:
\begin{equation}
t=m(q)=k_r=T_{36}\circ T_2\circ T_1\circ R(q).
\end{equation}
\begin{equation}
l=m\left(q^{1/5}\right)=k_{r/25}=T_{38}\circ T_2\circ T_4\circ T_1\circ R(q).
\end{equation}
Then
\begin{equation}
r=k_i(t)=k^{(-1)}(t)=\left(\frac{K\left(\sqrt{1-t^2}\right)}{K(t)}\right)^2.
\end{equation}
If
\begin{equation}
U^{*}(x)=\sqrt{\frac{-1}{2x^2}+\frac{x^4}{2}+\frac{\sqrt{1+18 x^6+x^{12}}}{2 x^2}}
\end{equation}
\begin{equation}
\alpha=U^{*}\left(\sqrt[12]{\frac{t\sqrt{1-t^2}}{l\sqrt{1-l^2}}}\right)^6
\end{equation}
\begin{equation}
p=2\sinh\left(\frac{1}{5}\textrm{arcsinh}\left(\frac{11+\alpha}{2}\right)\right)
\end{equation}
\begin{equation}
s=\sqrt[3]{\frac{(p-1)^5}{11+6p+6p^2+p^3+p^4}},
\end{equation}
the solution of
\begin{equation}
x^5-x-a=0\textrm{, }
\end{equation}
with
\begin{equation}
a=\frac{2(1+t^2)}{\sqrt[4]{5^5}\sqrt{t}\sqrt{1-t^2}}
\end{equation}
is
$$
x=4\left[\Psi(t,l)^{1/4}+l^{1/4}\right]Im\left[m\left(-(-q)^{1/5}\right)^{1/4}\zeta_{1/2}\right]\times
$$
\begin{equation}
\times Re\left[m\left((-1)^{2/5}q^{1/5}\right)^{1/4}\zeta_{1/4}\right],
\end{equation}  
$Y$ is the root of equation
\begin{equation}
Y^3+5s^{-1} Y^2-s Y-1=0
\end{equation}
and $\Psi(t,l)$ is given from Lemma 1 and Proposition 1.\\
\textbf{Remark.}\\
The interval of solvability $1<r\leq1+1/15$ can be increased if we use other branches of solutions of the equation
$$
256\frac{\left(y^2+(1-y^2)^2\right)^3}{y^4(1-y^2)^2}=x,
$$
which defines $T_{36},T_{38},T_{312}$.\\
\\
\textbf{Corollary.}\\
The Hermite's depressed equation is equivalent to Klein's equation and this is equivalent to solve the quintic equation.\\
\\
\textbf{Proof.}\\
Since $k_r=T_3(j(q))$ is defined from (23) and $T_1$ from (20) (first identity) and $T_2$ from Klein's equation (21) we can extract the result.

\subsection{Formulas for $T_1,T_{2},T_{36},T_{38},T_{312},T_4,T_5,T_6,T_7$ and $f_1$}

In this section we give closed form evaluation of algebraic functions $T_i,i=1,2,\ldots,7$.
$$
T_1(x):=-\frac{x^2}{3}+e^{-i\pi/3}\frac{3+x^5}{3 \left(18 x^4+x^9+3 \sqrt{3} \sqrt{-x^3+11 x^8+x^{13}}\right)^{1/3}}+
$$
\begin{equation}
+e^{i\pi/3}\frac{\left(18 x^4+x^9+3 \sqrt{3} \sqrt{-x^3+11x^8+x^{13}}\right)^{1/3}}{3 x}
\end{equation}
\\
\begin{equation}
T_2(x):=-\frac{\left(x^{20}-228 x^{15}+494 x^{10}+228 x^{5}+1\right)^{3}}{x^{5} \left(x^{10}+11 x^{5}-1\right)^{5}}
\end{equation} 
\\
\begin{equation}
T_{36}(x)=\sqrt{\frac{1}{2}-\frac{1}{2}\sqrt{1-4\left(\frac{768-x}{768}-e^{i\pi/3}\frac{x(1536-x)}{768D_1^{1/3}}+e^{-i\pi/3}\frac{D_1^{1/3}}{768}\right)}},
\end{equation}
\\
$$
T_{38}(x)=\sqrt{\frac{1}{2}+\frac{1}{2}\sqrt{1-4\left(\frac{768-x}{768}-e^{i\pi/3}\frac{x(1536-x)}{768D_1^{1/3}}+e^{-i\pi/3}\frac{Dz_1^{1/3}}{768}\right)}},\eqno{(85.1)}
$$
\\
$$
T_{312}(x)=\sqrt{\frac{1}{2}+\frac{1}{2}\sqrt{1-4\left(\frac{768-x}{768}-e^{-i\pi/3}\frac{x(1536-x)}{768D_1^{1/3}}+e^{i\pi/3}\frac{Dz_1^{1/3}}{768}\right)}},\eqno{(85.2)}
$$
where 
\begin{equation}
D_1=884736 x-2304 x^2+x^3+12288 \sqrt{3} \sqrt{1728 x^2-x^3}
\end{equation}
\\
\begin{equation}
T_4(x):=\left(\frac{x \left(1-2 x+4 x^2-3 x^3+x^4\right)}{1+3 x+4 x^2+2 x^3+x^4}\right)^{1/5}
\end{equation}
\\
\begin{equation}
T_5(x):=\frac{x}{\sqrt{x^2-1}}
\end{equation}
\\
\begin{equation}
T_6(x):=e^{i\pi/5}\frac{f_1\left(T_1(x)\right)}{xT_1(x)^2}
\end{equation}
where
$$
f_1(x):=-\frac{x^3}{3}-e^{-i\pi/3}\frac{\left(3x+x^6\right)}{3 \left(-18 x^4-x^9+3 \sqrt{3} \sqrt{-x^3+11 x^8+x^{13}}\right)^{1/3}}-
$$
\begin{equation}
-\frac{e^{i\pi/3}}{3} \left(-18 x^4-x^9+3 \sqrt{3} \sqrt{-x^3+11 x^8+x^{13}}\right)^{1/3}
\end{equation}
\\
$$
T_7(x):=-\frac{1}{3 x^2}-e^{i\pi/3}\frac{\left(-1+3 x^5\right)}{3 x^2 \left(1-18 x^5+3 \sqrt{3} \sqrt{-x^5+11 x^{10}+x^{15}}\right)^{1/3}}+
$$
\begin{equation}
+e^{-i\pi/3}\frac{\left(1-18 x^5+3 \sqrt{3} \sqrt{-x^5+11 x^{10}+x^{15}}\right)^{1/3}}{3 x^2}
\end{equation}
\\
$$
t_0(x)=\frac{1-\sqrt{1-x^2}}{1+\sqrt{1-x^2}}\eqno{(93.1)}
$$

\section{Appendix}

\subsection{The Tschirnhausen transform for the quintic}

For to reduce the general quintic 
$$
x^5+ax^4+bx^3+cx^2+dx+e=0\eqno{:(eq)}
$$
in the Bring quintic form we use the following Tschirnahausen transforms\\ 
\\
First use a quadratic transform
\begin{equation}
y=x^2+Ax+B
\end{equation}
with
\begin{equation}
A=\frac{4 a^4-13 a^2 b+15 a c\pm\sqrt{5\Delta_1}}{4 a^3-10 a b}
\end{equation}
and
\begin{equation}
B=\frac{5a^2b-20b^2+15a c\pm\sqrt{5\Delta_1}}{20a^2-50b},
\end{equation}
where
\begin{equation}
\Delta_1=-3 a^4 b^2+12 a^2 b^3+8 a^5 c-38 a^3 b c+45 a^2 c^2+16 a^4 d-40 a^2 b d
\end{equation}
and write $(eq)$ in the ''principal form'' 
\begin{equation}
x^5+A_1x^2+B_1x+C_1=0,
\end{equation}
where 
$$
A_1=r=3 a A b B-3 A^2 b B-3 b^2 B-6 a^2 B^2+6 a A B^2+12 b B^2-10 B^3-a A^2 c+A^3 c+
$$
\begin{equation}
+A b c+6 a B c-9 A B c-c^2-3 a A d+4 A^2 d+2 b d-6 B d-2 a e+5 A e
\end{equation}

$$
B_1=t=-3 a A b B^2+3 A^2 b B^2+3 b^2 B^2+4 a^2 B^3-4 a A B^3-8 b B^3+5 B^4+2 a A^2 B c-
$$
$$
-2 A^3 B c-2 A b B c-6 a B^2 c+9 A B^2 c+2 B c^2-a A^3 d+A^4 d+A^2 b d+6 a A B d-8 A^2 B d-
$$
\begin{equation}
-4 b B d+6 B^2 d-A c d+d^2-4 a A^2 e+5 A^3 e+3 A b e+4 a B e-10 A B e-2 c e
\end{equation}

$$
C_1=s=a A b B^3 - A^2 b B^3 - b^2 B^3 - a^2 B^4 + a A B^4 + 2 b B^4 - B^5 - 
 a A^2 B^2 c + A^3 B^2 c +
$$
$$
 + A b B^2 c + 2 a B^3 c - 3 A B^3 c - 
 B^2 c^2 + a A^3 B d - A^4 B d - A^2 b B d - 3 a A B^2 d + 
 4 A^2 B^2 d +
$$
$$
+ 2 b B^2 d - 2 B^3 d + A B c d - B d^2 - a A^4 e + 
 A^5 e + A^3 b e + 4 a A^2 B e - 5 A^3 B e -
$$
\begin{equation}
 - 3 A b B e - 2 a B^2 e + 
 5 A B^2 e - A^2 c e + 2 B c e + A d e - e^2
\end{equation}
Then again apply a Tschirnhausen transform (quartic)
\begin{equation}
y=x^4+kx^3+lx^2+mx+n,
\end{equation}
to reduce 
\begin{equation}
x^5+rx^2+tx+s=0
\end{equation}
into the ''Bring-Jerrard'' form 
\begin{equation}
x^5+A_2x+B_2=0
\end{equation}
Set 
\begin{equation}
\Delta_5=-27 r^6 s^2+256 r^2 s^5+108 r^7 t-1600 r^3 s^3 t+2250 r^4 s t^2+3125 r^2 t^4,
\end{equation} 
then we will have
\begin{equation}
k=\frac{-27 r^3 s+400 s^2 t-375 r t^2-3 \sqrt{5\Delta_5}}{54 r^4-320 s^3+600 r s t}
\end{equation}
\begin{equation}
l=\frac{18 r^3 s^2-45 r^4 t-250 r s t^2+2 s \sqrt{5\Delta_5}}{r \left(27 r^4-160 s^3+300 r s t\right)}
\end{equation}
\begin{equation}
n=\frac{135 r^4 s-1280 s^4+3600 r s^2 t-1125 r^2 t^2-9 r \sqrt{5\Delta_5}}{10 \left(27 r^4-160 s^3+300 r s t\right)}
\end{equation}
and $m$ is root of a 3rd degree equation
\begin{equation}
rX^3+\frac{M_1+N_1\sqrt{5\Delta_5}}{S_1}X^2+\frac{M_2+N_2\sqrt{5\Delta_5}}{S_2}X+\frac{M_3+N_3\sqrt{5\Delta_5}}{S_3}
\end{equation}
where
\begin{equation}
M_1=162 r^7-1104 r^3 s^3+2295 r^4 s t+1875 r^2 t^3,
\end{equation}

\begin{equation}
N_1=16 s^2-15 r t,
\end{equation}

\begin{equation}
S_1=54 r^5-320 r s^3+600 r^2 s t.
\end{equation}

$$
M_2=4374 r^{12}-59859 r^8 s^3+205440 r^4 s^6+124902 r^9 s t-869760 r^5 s^4 t+
$$
$$
+102400 r s^7 t+946350 r^6 s^2 t^2-688000 r^2 s^5 t^2+81000 r^7 t^3+510000 r^3 s^3 t^3+
$$
\begin{equation}
+1040625 r^4 s t^4+1250000 r s^2 t^5,
\end{equation}

\begin{equation}
N_2=-783 r^5 s^2+5120 r s^5+486 r^6 t-12960 r^2 s^3 t+4725 r^3 s t^2-10000 s^2 t^3
\end{equation}

\begin{equation}
S_2=2 r \left(27 r^4-160 s^3+300 r s t\right)^2.
\end{equation}

$$
M_3=-196830 r^{17}+4056885 r^{13} s^3-27812160 r^9 s^6+61649920 r^5 s^9+
$$
$$
+13107200 r s^{12}-8496495 r^{14} s t+117564615 r^{10} s^4 t-400708800 r^6 s^7 t-
$$
$$
-73728000 r^2 s^{10} t-126918900 r^{11} s^2 t^2+909441000 r^7 s^5 t^2-35200000 r^3 s^8 t^2-
$$
$$
-4829625 r^{12} t^3-651358125 r^8 s^3 t^3+426800000 r^4 s^6 t^3-133650000 r^9 s t^4+
$$
$$
+130500000 r^5 s^4 t^4+160000000 r s^7 t^4-1065234375 r^6 s^2 t^5+125000000 r^2 s^5 t^5-
$$
\begin{equation}
-10546875 r^7 t^6-828125000 r^3 s^3 t^6-263671875 r^4 s t^7,
\end{equation}

$$
N_3=50301 r^{10} s^2-675648 r^6 s^5+2298880 r^2 s^8-22599 r^{11} t+1658475 r^7 s^3 t-
$$
$$
-10376000 r^3 s^6 t-516375 r^8 s t^2+13329000 r^4 s^4 t^2-640000 s^7 t^2-
$$
$$
-2362500 r^5 s^2 t^3-4600000 r s^5 t^3-84375 r^6 t^4+
$$
\begin{equation}
+13375000 r^2 s^3 t^4-2109375 r^3 s t^5
\end{equation}

\begin{equation}
S_3=10 r (27 r^4 - 160 s^3 + 300 r s t)^3
\end{equation}
And
$$
A_2=5 n^4-2 m^3 n r+9 l m n^2 r-12 k n^3 r+2 l^3 n r^2-6 k l m n r^2+6 m^2 n r^2+9 k^2 n^2 r^2-
$$
$$
-9 l n^2 r^2-2 k^3 n r^3+6 k l n r^3-6 m n r^3+2 n r^4+5 l m^3 s-10 l^2 m n s-10 k m^2 n s+
$$
$$
+15 k l n^2 s+15 m n^2 s-2 l^4 r s+6 k l^2 m r s+3 k^2 m^2 r s-9 l m^2 r s-14 k^2 l n r s+
$$
$$
+16 l^2 n r s+2 k m n r s-12 n^2 r s+2 k^3 l r^2 s-6 k l^2 r^2 s+6 l m r^2 s+2 k n r^2 s-
$$
$$
-2 l r^3 s+5 k^2 l^2 s^2-5 l^3 s^2-5 k^3 m s^2-5 k l m s^2+5 m^2 s^2+10 k^2 n s^2+
$$
$$
+10 l n s^2-2 k^3 r s^2+4 k l r s^2-7 m r s^2+2 r^2 s^2+5 k s^3+m^4 t-8 l m^2 n t+
$$
$$
+6 l^2 n^2 t+12 k m n^2 t-16 n^3 t-l^3 m r t+3 k l m^2 r t-3 m^3 r t+2 k l^2 n r t-
$$
$$
-10 k^2 m n r t+4 l m n r t+15 k n^2 r t+k^3 m r^2 t-3 k l m r^2 t+3 m^2 r^2 t-
$$
$$
-2 k^2 n r^2 t+4 l n r^2 t-m r^3 t+k l^3 s t-7 k^2 l m s t+3 l^2 m s t+13 k m^2 s t+
$$
$$
+6 k^3 n s t-4 k l n s t-22 m n s t-3 k^4 r s t+11 k^2 l r s t-4 l^2 r s t-10 k m r s t+
$$
$$
+16 n r s t+3 k r^2 s t+k^2 s^2 t-6 l s^2 t+l^4 t^2-4 k l^2 m t^2+2 k^2 m^2 t^2+
$$
$$
+4 l m^2 t^2+8 k^2 l n t^2-8 l^2 n t^2-16 k m n t^2+18 n^2 t^2-k^3 l r t^2+3 k l^2 r t^2+
$$
$$
+k^2 m r t^2-5 l m r t^2-2 k n r t^2+l r^2 t^2+k^3 s t^2-3 k l s t^2+7 m s t^2-4 r s t^2+
$$
\begin{equation}
+k^4 t^3-4 k^2 l t^3+2 l^2 t^3+4 k m t^3-8 n t^3-k r t^3+t^4,
\end{equation}

$$
B_2=-n^5+m^3 n^2 r-3 l m n^3 r+3 k n^4 r-l^3 n^2 r^2+3 k l m n^2 r^2-3 m^2 n^2 r^2-
$$
$$
-3 k^2 n^3 r^2+3 l n^3 r^2+k^3 n^2 r^3-3 k l n^2 r^3+3 m n^2 r^3-n^2 r^4+m^5 s-
$$
$$
-5 l m^3 n s+5 l^2 m n^2 s+5 k m^2 n^2 s-5 k l n^3 s-5 m n^3 s-l^3 m^2 r s+
$$
$$
+3 k l m^3 r s-3 m^4 r s+2 l^4 n r s-6 k l^2 m n r s-3 k^2 m^2 n r s+9 l m^2 n r s+
$$
$$
+7 k^2 l n^2 r s-8 l^2 n^2 r s-k m n^2 r s+4 n^3 r s+k^3 m^2 r^2 s-3 k l m^2 r^2 s+
$$
$$
+3 m^3 r^2 s-2 k^3 l n r^2 s+6 k l^2 n r^2 s-6 l m n r^2 s-k n^2 r^2 s-m^2 r^3 s+
$$
$$
+2 l n r^3 s-l^5 s^2+5 k l^3 m s^2-5 k^2 l m^2 s^2-5 l^2 m^2 s^2+5 k m^3 s^2-
$$
$$
-5 k^2 l^2 n s^2+5 l^3 n s^2+5 k^3 m n s^2+5 k l m n s^2-5 m^2 n s^2-5 k^2 n^2 s^2-
$$
$$
-5 l n^2 s^2+k^3 l^2 r s^2-3 k l^3 r s^2-2 k^4 m r s^2+6 k^2 l m r s^2+3 l^2 m r s^2-
$$
$$
-7 k m^2 r s^2+2 k^3 n r s^2-4 k l n r s^2+7 m n r s^2-l^2 r^2 s^2+2 k m r^2 s^2-
$$
$$
-2 n r^2 s^2+k^5 s^3-5 k^3 l s^3+5 k l^2 s^3+5 k^2 m s^3-5 l m s^3-5 k n s^3-k^2 r s^3+
$$
$$
+2 l r s^3-s^4-m^4 n t+4 l m^2 n^2 t-2 l^2 n^3 t-4 k m n^3 t+4 n^4 t+l^3 m n r t-
$$
$$
-3 k l m^2 n r t+3 m^3 n r t-k l^2 n^2 r t+5 k^2 m n^2 r t-2 l m n^2 r t-5 k n^3 r t-
$$
$$
-k^3 m n r^2 t+3 k l m n r^2 t-3 m^2 n r^2 t+k^2 n^2 r^2 t-2 l n^2 r^2 t+m n r^3 t+
$$
$$
+l^4 m s t-4 k l^2 m^2 s t+2 k^2 m^3 s t+4 l m^3 s t-k l^3 n s t+7 k^2 l m n s t-
$$
$$
-3 l^2 m n s t-13 k m^2 n s t-3 k^3 n^2 s t+2 k l n^2 s t+11 m n^2 s t-k^3 l m r s t+
$$
$$
+3 k l^2 m r s t+k^2 m^2 r s t-5 l m^2 r s t+3 k^4 n r s t-11 k^2 l n r s t+
$$
$$
+4 l^2 n r s t+10 k m n r s t-8 n^2 r s t+l m r^2 s t-3 k n r^2 s t-k^4 l s^2 t+
$$
$$
+4 k^2 l^2 s^2 t-2 l^3 s^2 t+k^3 m s^2 t-7 k l m s^2 t+3 m^2 s^2 t-k^2 n s^2 t+
$$
$$
+6 l n s^2 t+k l r s^2 t-3 m r s^2 t+k s^3 t-l^4 n t^2+4 k l^2 m n t^2-2 k^2 m^2 n t^2-
$$
$$
-4 l m^2 n t^2-4 k^2 l n^2 t^2+4 l^2 n^2 t^2+8 k m n^2 t^2-6 n^3 t^2+k^3 l n r t^2-
$$
$$
-3 k l^2 n r t^2-k^2 m n r t^2+5 l m n r t^2+k n^2 r t^2-l n r^2 t^2+k^4 m s t^2-
$$
$$
-4 k^2 l m s t^2+2 l^2 m s t^2+4 k m^2 s t^2-k^3 n s t^2+3 k l n s t^2-7 m n s t^2-
$$
$$
-k m r s t^2+4 n r s t^2-l s^2 t^2-k^4 n t^3+4 k^2 l n t^3-2 l^2 n t^3-
$$
\begin{equation}
-4 k m n t^3+4 n^2 t^3+k n r t^3+m s t^3-n t^4
\end{equation}
\\
Lastly with a simple change of scaling we can reduce (102) further into the ''Bring quintic'' form
\begin{equation}
x^5+x+t=0
\end{equation}

\subsection{The Bring radicals hypergeometric function}

The Bring quintic (119) is solvable (not necessary in radicals), using the arguments of Lambert and Euler, (see [5] pg.306-307). Since it has the form
\begin{equation}
aqx^p+x^q=1,
\end{equation}
the next theorem solves (120) and hence (119) and (1) in terms of hypergeometric functions.\\
\\
\textbf{Theorem 3.2} (Lambert-Euler)\\The equation (120) admits root $x$ such that 
\begin{equation}
x^n=\frac{n}{q}\sum^{\infty}_{k=0}\frac{\Gamma(\{n+pk\}/q)(-q a)^k}{\Gamma(\{n+pk\}/q-k+1)k!}\textrm{, }n=1,2,3,\ldots
\end{equation}
where $\Gamma(x)$ is Euler's the Gamma function.\\

As a result of Theorem 3.2 we can define the hypergeometric function
\begin{equation}
\textrm{BR}(t)=-t\cdot{}_4F_3\left[\left\{\frac{1}{5},\frac{2}{5},\frac{3}{5},\frac{4}{5}\right\};\left\{\frac{1}{2},\frac{3}{4},\frac{5}{4}\right\};-\frac{3125 t^4}{256}\right].
\end{equation}
Then the solution of (119) will be $\textrm{BR}(t)$. Hence the solution of (102) is
\begin{equation}
x=\sqrt[4]{A_2}\cdot\textrm{BR}\left(\frac{B_2}{A_2^{5/4}}\right).
\end{equation}
For the solution $x$ of (102) also hold the next nested radicals expansion
\begin{equation}
x=\sqrt[5]{-B_2-A_2\sqrt[5]{-B_2-A_2\sqrt[5]{-B_2-A_2\sqrt[5]{-B_2-\ldots}}}}
\end{equation}

\subsection{The solvability of RRCF}

The same holds for the sextic equation (see [16]):
\begin{equation}
\frac{b^2}{20a}+by+ay^2=cy^{5/3} 
\end{equation}
Set 
\begin{equation}
Y=Y(q)=Y_r=R(q)^{-5}-11-R(q)^5,
\end{equation}
then (125) has solution
\begin{equation}
y=\frac{b}{250a}Y(q^2)=\textrm{, }q=e^{-\pi\sqrt{r}}\textrm{, }r>0,
\end{equation}
where $r$ can evaluated from the coefficients, using the relation $j_r=250\frac{c^3}{a^2b}$.\\ 
The $j$-invariant can be expressed in terms of the Ramanujan-Dedekind eta function as
\begin{equation}
j_r=\left[\left(q^{-1/24}\frac{\eta(\tau)}{\eta(2\tau)}\right)^{16}+16\left(q^{1/24}\frac{\eta(2\tau)}{\eta(\tau)}\right)^{8}\right]^3.
\end{equation}
The Ramanujan-Dedekind eta function is given from
\begin{equation}
\eta(\tau)=\prod^{\infty}_{n=1}\left(1-q^n\right)\textrm{, }q=e^{\pi i\tau}\textrm{, }Im(\tau)>0
\end{equation}
and special values can be found from the relation (see [9],[21])
\begin{equation}
\eta(\tau)^8=\frac{2^{8/3}}{\pi^4}q^{-1/3}(k_r)^{2/3}(k'_r)^{8/3}K(k_r)^4\textrm{, }q=e^{-\pi\sqrt{r}}\textrm{, }r>0
\end{equation} 
But most appropriate is to use (23) directly. 

\subsection{The Cubic Equation}

Suppose we have to solve the equation

\begin{equation}
x^3=3ax+b
\end{equation}
We assume that the solution is of the form 
\begin{equation}
x=\sqrt[3]{A}+\sqrt[3]{B}
\end{equation}  
Then 
$$
x^3=A+B+3\sqrt[3]{AB}\left(\sqrt[3]{A}+\sqrt[3]{B}\right)=A+B+3\sqrt[3]{AB}x
$$
Hence we have to solve 
\begin{equation}
A+B=b\textrm{ and }AB=a^3
\end{equation}
which reduces clearly to quadratic equation
\begin{equation}
X^2-bX+a^3=0
\end{equation}
\[
\]

\centerline{\bf References}\vskip .2in

[1]: I.J.Zucker. 'The summation of series of hyperbolic functions'. SIAM J. Math. Ana.10.192, 1979.

[2]: G.E.Andrews. 'Number Theory'. Dover Publications, New York, 1994.

[3]: T.Apostol. 'Introduction to Analytic Number Theory. Springer Verlang, New York, Berlin, Heidelberg, Tokyo, 1974.

[4]: M.Abramowitz and I.A.Stegun. 'Handbook of Mathematical Functions'. Dover Publications, 1972.

[5]: Bruce.C.Berndt. 'Ramanujan`s Notebooks Part I'. Springer Verlag, New York, 1985.

[6]: Bruce.C.Berndt. 'Ramanujan`s Notebooks Part II'. Springer Verlag, New York, 1989.

[7]: Bruce.C.Berndt. 'Ramanujan's Notebooks Part III'. Springer Verlag, New York, 1991.

[8]: Bruce.C.Berndt, Sen-Shan Huang, Jaebum Sohn and Seung Hwan Son. 'Some Theorems on the Rogers-Ramanujan Continued Fraction in Ramanujan Lost Notebook'. Trans. Amer. Math. Soc. 352 (2000), 2157-2177

[9]: E.T.Whittaker and G.N.Watson. 'A course on Modern Analysis'. Cambridge U.P., 1927.

[10]: Nikos Bagis. 'The complete evaluation of Rogers Ramanujan and other continued fractions with elliptic functions'. arXiv:1008.1304v1 [math.GM], 2010.

[11]: Michael Trott. 'Modular Equations of the Rogers-Ramanujan Continued Fraction'. page stored in the Web.

[12]: Nikos Bagis. 'The First Derivative of Ramanujan Cubic Fraction'. arXiv:1103.5346v1 [math.GM], 2011.  

[13]: Nikos Bagis. 'Parametric Evaluations of the Rogers-Ramanujan Continued Fraction'. IJMMS, Vol. 2011.

[14]: Nikos Bagis. 'The $w$-modular function and the evaluation of Rogers-Ramanujan continued fraction'. International Journal of Pure and Applied Mathematics. Vol 84, No 1, 2013, 159-169.

[15]: J.M.Borwein and P.B.Borwein. 'Pi and the AGM'. John Wiley and Sons, Inc. New York, Chichester, Brisbane, Toronto, Singapore, 1987.

[16]: Nikos Bagis. 'On a General Sextic Equation Solved by the Rogers-Ramanujan Continued Fraction'. arXiv:1111.6023v2 [math.GM], 2012. 

[17]: Soon-Yi Kang. 'Ramanujan's formulas for the explicit evaluation of the Rogers-Ramanujan continued fraction and theta functions'. Acta Arithmetica. XC.1, (1999).

[18]: Bruce.C.Berndt. 'Ramanujan Notebooks Part V'. Springer Verlag, New York (1998).

[19]: Bruce.C.Berndt, Heng Huat Chan, Sen-Shan Huang, Soon-Yi Kang, Jaebum Sohn and Seung Hwan Son. 'The Rogers-Ramanujan Continued Fraction'. J. Comput. Appl. Math., 105 (1999), 9-24.

[20]: J.V.Armitage, W.F.Eberlein. 'Elliptic Functions'. Cambridge University Press. (2006)

[21]: J.M. Borwein, M.L. Glasser, R.C. McPhedran, J.G. Wan, I.J. Zucker. 'Lattice Sums Then and Now'. Cambridge University Press. New York, (2013).

[22]: Heinrich Dorrie. '100 Great Problems of Elementary Mathematics (Their History and Solution)'. Dover Publications, inc. New York. (1965) 

[23]: N.N. Lebedev. 'Special Functions and their Applications'. Dover Pub. New York. (1972)

[24]: W. Duke. 'Continued fractions and Modular functions'. Bull. Amer. Math. Soc. (N.S.), 42 (2005), 137-162.

[25]: Nikos Bagis. 'Evaluation of Fifth Degree Elliptic Singular Moduli'. arXiv:1202.6246v1. (2012)

[26]: Nikos Bagis. 'Generalizations of Ramanujan's Continued Fractions'.\\ arXiv:11072393v2 [math.GM] 7 Aug 2012.

\end{document}